\documentclass[12pt]{amsart}
\input{diagrams}

\newtheorem{proposition}{Proposition}[section]
\newtheorem{lemma}[proposition]{Lemma}

\newtheorem{theorem}[proposition]{Theorem}

\theoremstyle{definition}
\newtheorem{definition}[proposition]{Definition}

\theoremstyle{remark}

\newcommand{\thlabel}[1]{\label{th:#1}}
\newcommand{\thref}[1]{Theorem~\ref{th:#1}}
\newcommand{\selabel}[1]{\label{se:#1}}
\newcommand{\seref}[1]{Section~\ref{se:#1}}
\newcommand{\lelabel}[1]{\label{le:#1}}
\newcommand{\leref}[1]{Lemma~\ref{le:#1}}
\newcommand{\prlabel}[1]{\label{pr:#1}}
\newcommand{\prref}[1]{Proposition~\ref{pr:#1}}

\newcommand{\delabel}[1]{\label{de:#1}}
\newcommand{\deref}[1]{Definition~\ref{de:#1}}
\newcommand{\eqlabel}[1]{\label{eq:#1}}
\newcommand{\equref}[1]{(\ref{eq:#1})}

\newcommand{\Hom}{{\rm Hom}}
\newcommand{\End}{{\rm End}}

\def\ot{\otimes}

\newcommand{\Mm}{\mathcal{M}}

\def\text#1{{\rm {\rm #1}}}

\def\smashco{\mathrel>\joinrel\mathrel\triangleleft}

\def\os{\otimes}

\def\De{\Delta}
\def\si{\sigma}
\def\ga{\gamma}
\def\a{\alpha}

\def\ve{\varepsilon}
\def\lra{\longrightarrow}

\def\cd{\cdot}

\def\la{\lambda}

\def\ds{\diamondsuit}

\def\cbc{C\Box _B C}
\def\rh{{\rm Hom}}

\def\crth{C\smashco _{\a'}' H}

\def\ct{{\mathcal T}}
\def\cl{{\mathcal L}}
\def\bs{{\bar S}}

\begin{document}
\title[Twistings and Hopf-Galois coextensions]
{Twistings, crossed coproducts and Hopf-Galois coextensions}
\author{S. Caenepeel}
\address{Faculty of Applied Sciences,
Vrije Universiteit Brussel, VUB, B-1050 Brussels, Belgium}
\email{scaenepe@vub.ac.be}
\urladdr{http://homepages.vub.ac.be/~scaenepe/}
\author{Dingguo Wang}
\address{Department of Mathematics,
Qufu Normal University, Qufu, Shandong 273165, China}
\email{diwang@vub.ac.be}
\author{Yanxin Wang}
\address{Department of Mathematics,
Tsinghua University, Beijing 100084, China}
\thanks{Research supported by the project G.0278.01 ``Construction
and applications of non-commutative geometry: from algebra to physics"
from FWO Vlaanderen}
\urladdr{}
\begin{abstract}
Let $H$ be a Hopf algebra. Ju and Cai
introduced the notion of twisting of an $H$-module coalgebra.
In this note, we study the relationship between twistings,
crossed coproducts and Hopf-Galois coextensions. In particular,
we show that a twisting of an $H$-Galois coextension
remains $H$-Galois if the twisting is invertible.
\end{abstract}
\maketitle

\section*{Introduction}\selabel{0}
A fundamental result in Hopf-Galois theory is the Normal Basis
Theorem, stating that, for a finitely generated cocommutative
Hopf algebra $H$ over a commutative ring $k$, the set of isomorphism
classes of Galois $H$-objects that are isomorphic to $H$
as an $H$-comodule is a group, and this group is isomorphic
to the second Sweedler cohomology group $H^2(H,k)$
(see \cite{Sw}). The Galois
object corresponding to a 2-cocycle is then given by a
crossed product construction. The crossed product construction
can be generalized to arbitrary Hopf algebras, and
plays a fundamental role in the theory of extensions of
Hopf algebras, see \cite{BCM} and \cite{DT}. Also in this
more general situation, it turns out that there is a close
relationship between crossed products on one side, and
Hopf-Galois extensions and cleft extensions, cf.
\cite{BCM}, \cite{BCM}, \cite{DT}. A survey can be found in
\cite{Mon}. An alternative way to deform
the multiplication on an $H$-comodule algebra $A$ has been proposed
in \cite{BCZ}, using a so-called twisting of $A$, and it was
shown that the crossed product construction can be viewed as
a special case of the twisting construction. The relation
between twistings and $H$-Galois extensions was studied in
\cite{BT}.\\
Now there exists a coalgebra version of the Normal Basis
Theorem (see \cite{C}). In this situation, one tries to deform the
comultiplicationon a commutative Hopf algebra $H$, using this time a
Harrison cocycle instead of a Sweedler cocycle. 
Crossed coproducts, cleft coextensions and
Hopf-Galois coextensions have been introduced and studied in
\cite{DMR} and \cite{DRZ}. Ju and Cai \cite{JC} have 
introduced the notion twisting of an $H$-module coalgebra,
which can be viewed as dual version of the twistings
introduced in \cite{BCZ}. The aim of this paper is to 
study the relationship between twistings, crossed coproducts
and Hopf-Galois coextensions. Our main result is the fact
that the twisting of a Hopf-Galois coextension by an invertible
twist map is again a Hopf-Galois coextension (and conversely).\\
Our paper is set up as follows: in \seref{1.1}, we recall
the twistings introduced in \cite{JC}, and in \seref{1.2}
the definition of a Harrison cocycle and
the crossed coproduct construction from \cite{DMR} and \cite{DRZ}.
In \seref{2}, we introduce an alternative version of 2-cocycles,
called twisted 2-cocycles, and discuss the relation with
Harrison cocycles (\prref{3.4}). In \seref{3}, we introduce an equivalence
relation on the set of twistings of an $H$-module coalgebra,
and we show that a twisting in an equivalence class is invertible
if and only if all the other twistings in this equivalence
class are invertible  (\thref{2.5}). Two twistings are
equivalentif and only if their corresponding crossed coproducts are
isomorphic (\prref{1.2}). In \seref{4} the relationship between
twistings and Hopf-Galois coextensions is investigated.\\
For the general theory of Hopf algebras, we refer to the literature,
see for example \cite{DNR}, \cite{Mon}, \cite{Swe}.

\section{Notation and preliminary results}\selabel{1}
We work over a field $k$. All maps are assumed to be $k$-linear.
For the comultiplication on a $k$-coalgebra $C$, we use the
Sweedler-Heyneman notation
$$\Delta_C(c)=c_1\ot c_2$$
with the summation implicitely understood. We use a similar notation
for a (right) coaction of a coalgebra on a comodule:
$$\rho(m)=m_0\ot m_1\in M\ot C$$
Let $A$ be a $k$-algebra,
then $\Hom(C,A)$ is also an algebra, with convolution product
$$(f*g)(c)=f(c_1)g(c_2)$$
${\rm Reg}(C,A)$ will denote the set of convolution invertible
elements in $\Hom(C,A)$. ${\mathcal M}_A^C$ will be the category
of modules with a right $A$-action and a right $C$-coaction,
such that the $C$-coaction is $A$-linear.

\subsection{Twistings of a coalgebra}\selabel{1.1}
We recall some definitions and results from \cite{JC}.
Let $H$ be a Hopf algebra over a field $k$, with bijective antipode $S$.
The composition inverse of the antipode will be denoted by $\bs$.\\
Recall that a right $H$-module coalgebra is a
coalgebra $C$ which is also a right $H$-module such that
$$\Delta(c\cdot h)=c_1\cdot h_1\ot c_2\cdot h_2~~{\rm and}~~
\varepsilon_C(c\cdot h)=\varepsilon_C(c)\varepsilon_H(h)$$
for all $c\in C$ and $h\in H$.\\
${\mathcal M}_H^C$ is the category whose objects are right $H$-modules
and right $C$-comodules $M$ such that the following compatibility
relation is satisfied:
$$\rho(m\cdot h)=m_0\cdot h_1\ot m_1\cdot h_2$$
Recall from \cite{JC} that we have the following associative
multiplication on $\Hom(C, H\ot C)$:
$$\tau\ast\la =(m_H\os id_C)\circ(id_H\os \la)\circ\tau$$
for all $\tau,\lambda\in \Hom(C, H\ot C)$. The unit of this multiplication
is the map $\sigma:\ C\to H\ot C$, $\sigma(c)=1\ot c$.\\
Remark that we have an algebra isomorphism
$$\alpha:\ \Hom(C, H\ot C)\to {}_H\End(H\ot C)^{\rm op}$$
For $\tau:\ C\to H\ot C$, we define the corresponding
$\alpha(\tau)=f_{\tau}:\ H\ot C\to H\ot C$ by
$$f_{\tau}(h\ot c)=h\tau(c)=hc_{-1}\ot c_{0}$$
Assume that $\tau$ satisfies the following normality conditions:
\begin{equation}\eqlabel{4}
(1\os \ve_C)\tau(c)=\ve(c)1_H,~~~ (\ve_H\os 1)\tau(c)=c
\end{equation}
If we write $\tau(c)= c_{-1}\os c_0$ (summation understood), then \equref{4} takes
the following form
$$ c_{-1}\varepsilon_C(c_0)=\varepsilon_C(c)1_H,~~~
\varepsilon_H(c_{-1})c_0=c$$
We can then define a new (in general non-coassociative)
comultiplication $\De_\tau$ on $C$ as follows:
$$\De_\tau(c)= c_1\cd c_{2.-1}\os c_{2.0},~~{\rm or}~~
\Delta_{\tau}=(m_H\ot id)\circ (id\ot \tau)\circ \Delta$$
Let $C^{\tau}$ be equal to $C$ as a right $H$-module, with comultiplication
$\De_\tau(c)$. A similar construction applies to $M\in {\mathcal M}_H^C$:
$F_{\tau}(M)=M^{\tau}$ as a right $H$-comodule, with 
$$\rho^{\tau}(m)=m_0\cd m_{1.-1}\os m_{1.0}$$
$\tau$ is called a twisting if and only if $C^\tau$ is a right $H$-module
coalgebra, and $M^{\tau}\in\Mm^{C^\tau}_H$ for all $M\in \Mm^{C}_H$.
It is shown in \cite[Theorem 1.1]{JC} that 
$\tau:\ C\to H\ot C$ satisfying \equref{4} is a twisting
if and only if for all $h\in H$ and $c\in C$,
\begin{equation}
 c_{-1}h_1\os c_0\cd h_2= h_1(c\cd h_2)_{-1}\os (c\cd
h_2)_0\eqlabel{6}
\end{equation}
and
\begin{equation}
 c_{-1}\os c_{0.1}\cd c_{0.2.-1}\os
c_{0.2.0}= c_{1.-1}c_{2.-1.1}\os c_{1.0}\cd c_{2.-1.2}
\os c_{2.0}\eqlabel{7}
\end{equation}
\equref{6} is equivalent to
\begin{equation} \eqlabel{10}
 S(h_1)c_{-1}h_2\os c_0\cd
h_3= (c\cd h)_{-1}\os (c\cd h)_0
\end{equation}
If $\tau$ has an inverse $\lambda$, then the functor $F$ is an
equivalence of categories.\\
Left hand twistings are defined in a similar way. Consider the vector space isomorphism
$$\Hom(C,C\ot H)\cong \End_H(C\ot H^{\rm op},C\ot H^{\rm op})$$
The composition on the right hand side is transported into the
following associative multiplication on $\Hom(C,C\ot H)$:
$$ \tau\times \la =T\circ(T\circ\la \ast T\circ \tau)$$
Here $T$ is the usual twist map. The unit $\sigma'$ on
$\Hom(C,C\ot H)$ is given by $\sigma'(c)=c\ot 1$.
If $\lambda\in \Hom(C,C\ot H)$
satisfies the normalizing conditions
\begin{equation}\eqlabel{5}
(1\os \ve_H)\la(c)
 =c,~~~ (\ve_C\os 1)\la(c)=\ve_C(c)1_H
\end{equation}
then we can twist the comultiplication on $C$ as follows:
write $\la(c)= c_0\os c_1$, and define $ {}_\la \De$ by
$${}_\la\De(c)= c_{1.0}\os c_2\cd c_{1.1}$$
${}^{\lambda}C$ will be $C$ as a right $H$-module, with the
comultiplication $ {}_\la \De$. The $C$-coaction $M\in {}^C{\mathcal M}_H$
can also be twisted:
$${}^{\lambda}\rho(m)=m_{-1.0}\ot m_0m_{-1.1}$$
$\lambda$ is called a left hand twisting if ${}^{\lambda}C$ is
an $H$-module coalgebra, and ${}^{\lambda}M\in {}^{{}^{\lambda}C}{\mathcal M}_H$
for every $M\in {}^C{\mathcal M}_H$.
$\lambda:\ C\to C\ot H$ satisfying \equref{5} is a 
left hand twisting
if and only if for all $h\in H$ and $c\in C$,
\begin{equation}
 c_0\cd h_1\os c_1h_2= (c\cd h_1)_0\os h_2(c\cd
h_1)_1\eqlabel{8}
\end{equation}
and
\begin{equation}
 c_{0.1.0}\os c_{0.2}\cd c_{0.1.1}\os c_1= c_{1.0}\os c_{2.0}\cd c_{1.1.1}\os c_{2.1}c_{1.1.2}\eqlabel{9}
\end{equation}
\equref{8} is equivalent to 
\begin{equation} \eqlabel{11} 
 c_0\cd h_1\os \bar{S}(h_3)c_1h_2=\sum(c\cd h)_0\os (c\cd h)_1
\end{equation} 

For $\tau \in \Hom(C, H\ot C)$ with inverse $\lambda$, we write
\begin{equation}\eqlabel{notation}
\tau(c)= c_{-1}\os c_0,~~~
\la(c)= c_{(-1)}\os c_{(0)}
\end{equation}
We then have 
\begin{equation}\eqlabel{12} 
 c_{-1}c_{0.(-1)}\os c_{0.(0)}=c_{(-1)}c_{(0).-1}\os c_{(0).0}
=1\os c
\end{equation}
For $\gamma \in \Hom(C, C\ot H)$ with inverse $\mu$, we write
$$\ga(c)= c_0\os c_1,~~~\mu(c)= c_{(0)}\os c_{(1)}$$
Let ${\mathcal T}(C)$ and ${\mathcal L}(C)$ be the sets of respectively
twistings and left hand twistings of $C$, and 
$U(\ct(C))$, $U(\cl(C))$ the sets of invertible 
twistings and left hand twistings.

\begin{proposition}\prlabel{1.2}
Take $\tau\in U(\ct(C))$ with inverse $\lambda$. Define
$\ell(\tau):\ C\lra C\os H$ by 
$$\ell(\tau)(c)= c_{0.(0)}\cd
\bs(c_{0.(-1)})\bs(c_{-1})_1 \os \bs(c_{-1})_2$$
Take $\ga\in U(\cl(C))$, with inverse $\mu$. Define
$r(\ga):\ C\lra H \os C$ by 
$$r(\ga)(c)= S(c_1)_1\os c_{0.(0)}\cd S(c_{0.(1)}) S(c_1)_2$$
Then $\ell:\ U(\ct(C))\to U(\cl(C))$ is a bijection with
inverse $r$. Furthermore 
$\ell(\si)=\si'$ and $r(\si')=\si$.
\end{proposition}

\begin{proof}
It is shown in \cite{JC} that $\ell(\tau)\in U(\cl(C))$ with
inverse given by
$$ \ell(\tau)'(c) = c_{0.(0)}\cd
\bs(c_{0.(-1)})_1\bs(c_{-1})_1 \os
\bs(\bs(c_{-1})_3)\bs(c_{0.(-1)})_2 \bs(c_{-1})_2$$
Set
$g=\bs(c_{0.(-1)}),\ h=\bs(c_{-1})$. Then $\ell(\tau)(c)= c_{0.(0)}\cd gh_1\os h_2$, so
\begin{eqnarray*}
&&\hspace*{-2cm}
 r(\ell(\tau))(c)= S(h_2)_1\os (c_{0.(0)}\cd gh_1)_{0.(0)}
 \cd \bs((c_{0.(0)}\cd gh_1)_{0.(-1)})_1\\
&&\hspace*{5mm}
\bs((c_{0.(0)}\cd gh_1 )_{-1})_1 
S(\bs(\bs((c_{0.(0)}\cd gh_1)_{-1})_3)\\
&&\hspace*{5mm}
\bs((c_{0.(0)}\cd gh_1)_{0.(-1)})_2
\bs((c_{0.(0)}\cd gh_1)_{-1})_2)   S(h_2)_2\\
&=& S(h_2)_1\os (c_{0.(0)}\cd gh_1)_{0.(0)}
 \bs((c_{0.(0)}\cd gh_1)_{0.(-1)})_1\\
&&\hspace*{5mm}
\bs((c_{0.(0)}\cd gh_1
 )_{-1})_1
S(\bs((c_{0.(0)}\cd gh_1)_{-1})_2)\\
&&\hspace*{5mm}
S(\bs((c_{0.(0)}\cd gh_1)_{0.(-1)})_2)
\bs((c_{0.(0)}\cd gh_1)_{-1})_3S(h_2)_2\\
&=& S(h_2)_1\os (c_{0.(0)}\cd gh_1)_{0.(0)}\\
&&\hspace*{5mm}
 \cd \ve((c_{0.(0)}\cd gh_1)_{0.(-1)}\bs((c_{0.(0)}\cd
 gh_1)_{-1})S(h_2)_2\\
&=& S(h_2)_1\os (c_{0.(0)}\cd gh_1)_0
 \cd \bs((c_{0.(0)}\cd gh_1)_{-1})S(h_2)_2\\
&=& S(h_2)_1\os c_{0.(0).0}\cd (gh_1)_3\bs((gh_1)_2)
\bs(c_{0.(0).-1})\\
&&\hspace*{5mm}  \bs(S(gh_1)_1)S(h_2)_2\\ 
&=& S(h_3)\os
c_{0.(0).0}\cd \bs(c_{0.(0).-1})gh_1S(h_2)\\ 
&=& S(h)\os
c_{0.(0).0}\cd \bs(c_{0.(0).-1})g
\end{eqnarray*}
\begin{eqnarray*} 
&=& S(\bs(c_{-1}))\os
c_{0.(0).0}\cd \bs(c_{0.(-1)}c_{0.(0).-1}) \\
&=& c_{-1}\os
c_0=\tau(c)
\end{eqnarray*}
In \cite{JC}, it is also shown that $r(\ga) \in U(\ct(C))$, and it
is straightforward to verify that the $\ast$-inverse of
$r(\ga)$ is given by
$$r(\ga)'(c)= S(S(c_1)_1)S(c_{0.(1)})_1S(c_1)_2\os c_{0.(0)}\cd
S(c_{0.(1)})_2S(c_1)_3$$
A routine verification similar to the one above then shows that
$$\ell(r(\ga))(c)=\gamma(c)$$
for all $c\in C$.
It is easy to show that $\ell(\si)=\si'$ and $r(\si')=\si$.
\end{proof}

\subsection{The crossed coproduct}\selabel{1.2}
We recall the following definitions from \cite{DMR} and \cite{DRZ}.

\begin{definition}\delabel{1.3}
Let $C$ be a coalgebra and $H$ a Hopf algebra.
We say that $H$ coacts weakly on $C$ if there is a
$k$-linear map $\rho:\ C\lra H\os C;\ \rho(c)= c_{[-1]}\os
c_{[0]}$ satisfying the following conditions, for all $c\in C$: 
\begin{eqnarray}
 c_{[-1]}\os c_{[0]1}\os c_{[0]2}&=&
 c_{1[-1]}c_{2[-1]}\os c_{1[0]} \os c_{2[0]}\eqlabel{14}\\ 
 \ve_C(c_{[0]})c_{[-1]}&=&\ve(c)1_H\eqlabel{15} \\ 
 \ve_H(c_{[-1]})c_{[0]}&=&c\eqlabel{16}
\end{eqnarray}
\end{definition}

Assume that $H$ coacts weakly on $C$, and let $\a:\ C\lra H\os H;\
\a(c)= \a_1(c)\os \a_2(c)$ be a linear map. Let
$C\smashco_{\alpha} H$ be the coalgebra whose
underlying vector space is $C\os H$, with
comultiplication and counit given by 
\begin{eqnarray*}
\De_\a(c\smashco h)&=& (c_1\smashco
c_{2[-1]}\a_1(c_3)h_1)\os
(c_{2[0]}\smashco \a_2(c_3)h_2)\\ 
\ve_\a(c\smashco h)&=&\ve_C(c)\ve_H(h)
\end{eqnarray*}
It was pointed out in \cite{DRZ} that $\ve_\a(c\smashco h)$
satisfies the counit property if and only if
\begin{equation}\eqlabel{I}
(\varepsilon_H\ot id)\alpha(c)=(id\ot \varepsilon_H)\alpha(c)=\varepsilon_C(c)1_H
\end{equation}
$\De_\a$ is coassociative if and only if $\alpha$ satisfies
\begin{eqnarray}
&&\hspace*{-15mm}  c_{1[-1]}\a_1(c_2)\os
\a_1(c_{1[0]})\a_2(c_2)_1\os
\a_2(c_{1[0]})\a_2(c_2)_2\nonumber\\
&=&  \a_1(c_1)\a_1(c_2)_1\os
\a_2(c_1)\a_1(c_2)_2\os \a_2(c_2)\eqlabel{II}\\
&&\hspace*{-15mm}  c_{1[-1]}\a_1(c_2)\os c_{1[0][-1]}\a_2(c_2)\os
c_{1[0][0]}\nonumber\\ 
&=& \a_1(c_1)c_{2[-1]1}\os
\a_2(c_1)c_{2[-1]2} \os c_{2[0]}\eqlabel{III}
\end{eqnarray}
In \cite{DRZ}, \equref{II} is called the cocycle
condition, and \equref{III} is called the twisted
comodule condition. Following \cite{CDMP}, we call
$\alpha$ satisfying (\ref{eq:I}-\ref{eq:III}),
a Harrison 2-cocycle.\\

Now consider two weak $H$-coactions $\rho,\ \rho':\ C\lra H\os C$,
and write
$$\rho(c)= c_{[-1]}\os c_{[0]}~~{\rm and}~~\rho'(c)=  c_{<-1>}\os c_{<0>}$$
Also consider two 2-cocycles $\alpha,\alpha':\ C\lra H\os H$
corresponding respectively to $\rho$ and $\rho'$, and write
$$\a(c)= \a_1(c)\os \a_2(c)~~{\rm and}~~\a'(c)= \a'_1(c)\os \a'_2(c)$$
Then we can consider the crossed coproducts
$C\smashco _{\alpha} H$ and $C\smashco _{\a'}^{'} H$. In the next
Lemma, we discuss when these are isomorphic.

\begin{lemma}\lelabel{1.4}
Consider a convolution invertible map $u:\ C\to H$
satisfying the conditions
\begin{eqnarray}
&&  c_{<-1>}\os
c_{<0>}= u^{-1}(c_1)c_{2[-1]}u(c_3)\os c_{2[0]}\eqlabel{18}\\ 
&&\a'(c)= u^{-1}(c_1)c_{2[-1]}\a_1(c_3)u(c_4)_1\os
u^{-1}(c_{2[0]})\a_2(c_{3})u(c_4)_2\eqlabel{19}
\end{eqnarray}
for all $c\in C$. Then the map
\begin{equation}
\phi:\ \crth\lra C\smashco_{\alpha} H;~~
\phi(c\smashco' h)= c_1\smashco u(c_2)h\eqlabel{17}
\end{equation}
is a left $C$-colinear, right $H$-linear coalgebra isomorphism.
Every left $C$-colinear, right $H$-linear coalgebra isomorphism
between $C\smashco _{\alpha} H$ and $C\smashco _{\a'}^{'} H$
is of this type.
\end{lemma}

\begin{proof}
The proof is a dual version of a similar statement for crossed products,
see \cite{Mon}.
\end{proof}

It was shown in \cite{JC} that
the crossed coproduct construction can be viewed as a special case
of the twisting construction from \seref{1.1}.
Let $H$ be a Hopf algebra, and $C$ a right $H$-module coalgebra,
and view
$C\otimes H$ as a right $H$-module coalgebra, with the right $H$-action
is induced by the multiplication in $H$. It was proved in \cite{JC}
that there is a bijective correspondence between crossed coproduct
structures on $C\otimes H$ and twistings of $C\otimes H$. Let us
recall the description of this bijection.\\
Consider a weak coaction $\rho$ and a 2-cocycle $\alpha$ giving
rise to the crossed coproduct $C\smashco _{\alpha} H$, and write
$$\rho(c)= c_{[-1]}\os c_{[0]}~~;~~\a(c)= \a_1(c)\os \a_2(c)$$
The corresponding twisting $\tau:\ C\os H\lra H\os C\os H$ is defined by
\begin{equation}\eqlabel{27}
\tau(c\os h)= S(h_1)c_{1[-1]}\a_1(c_2)h_2\os c_{1[0]}\os \a_2(c_2)h_3
\end{equation}
Conversely, if $\tau$ is a twisting of $C\os H$, then $(C\os
H)^\tau=C\smashco_\alpha H$, with weak coaction $\rho$
and 2-cocycle $\alpha$ given by 
\begin{eqnarray} 
\rho(c)&=&(id\os id\os \ve_H)\tau(c\os 1)\eqlabel{28}\\
\a(c)&=&(id\os\ve_C\os id)\tau(c\os 1)\eqlabel{29}
\end{eqnarray}

\section{Twisted 2-cocycles}\selabel{2}
Let $H$ be a Hopf algebra with bijective antipode $S$; let
$\bs$ be the composition inverse of $S$. Take an $H$-module
coalgebra $C$, and let $B=C/CH^+$.

\begin{definition}\delabel{2.1}
A map $\a:\ C\lra H\os H,~~ \a(c)= \a_1(c)\os \a_2(c)$ is called a twisted 2-cocycle if 
the following conditions are satisfied, for all $h\in H$ and $c\in
C$: 
\begin{eqnarray}
&&(id_H\os \ve_H)\a(c)=(\ve_H\os id_H)\a(c)=\ve_C(c)1_H\eqlabel{a}\\
&& \a(c\cd h)= S(h_1)\a_1(c)h_2\os
\bs(h_4)\a_2(c)h_3\eqlabel{b}\\
\eqlabel{c}&& \a_1(c_1)\a_1(c_3)_1\os c_2\cd
\a_2(c_1)\a_1(c_3)_2\os
\a_2(c_3)\\ 
&&\hspace*{1cm}= \a_1(c_1)\os c_2\cd \a_1(c_3)\a_2(c_1)_1\os
\a_2(c_3)\a_2(c_1)_2\nonumber
\end{eqnarray}
\end{definition}

Our first result is the fact that twisted 2-cocycles can be used to
define twistings on $C$.

\begin{proposition}\prlabel{2.2}
With notation as above, if $\alpha:\ C\to H\ot H$ is a twisted
2-cocycle, then the map 
$$\tau_\a:\ C\lra H\os C;~~\tau_\a(c)=  \a_1(c_1)\os c_2\cd \a_2(c_1)$$
is a twisting of $C$.
\end{proposition}

\begin{proof}
It follows easily from \equref{a} that $\tau_{\alpha}$ satisfies
the normalizing condition \equref{4}. Next we compute that
\begin{eqnarray*} 
&&\hspace*{-2cm}
 (c\cd h)_{-1}\os (c\cd h)_0
= \a_1((c\cd h)_1)\os
(c\cd h)_2\cd \a_2((c\cd h_2)_1) \\ 
&=& \a_1(c_1\cd h_1)\os
c_2\cd h_2\a_2(c_1\cd h_1)\\ 
&=& S(h_1)\a_1(c_1)h_2\os c_2\cd
h_5\bs(h_4)\a_2(c_1)h_3\\ 
&=& S(h_1)\a_1(c_1)h_2\os c_2\cd
\a_2(c_1)h_3\\ 
&=& S(h_1)c_{-1}h_2\os c_0\cd h_3
\end{eqnarray*}
and \equref{6} follows easily. Finally we compute the left and right
hand side of \equref{7}
\begin{eqnarray*}
&&\hspace*{-2cm}
 c_{-1}\os c_{0.1}\cd c_{0.2.-1}\os
c_{0.2.0}= 
(1\os \De_\tau)\tau_\a(c)\\
&=& \a_1(c_1)\os
(c_2\cd \a_2(c_1))_1 \cd \a_1(((c_2\cd \a_2(c_1))_2)_1)\\
&&\hspace*{5mm}
\os
((c_2\cd \a_2(c_1))_2)_2\cd \a_2(((c_2\cd \a_2(c_1))_2)_1)\\ 
&=& \a_1(c_1)\os c_2\cd \a_2(c_1)_1\a_1(c_3\cd \a_2(c_1)_2)\\
&&\hspace*{5mm} 
\os(c_4\cd
\a_2(c_1)_3)\cd \a_2(c_3\cd \a_2(c_1)_2)\\ 
&=& \a_1(c_1)\os
c_2\cd \a_2(c_1)_1S(\a_2(c_1)_2)\a_1(c_3) \a_2(c_1)_3\\
&&\hspace*{5mm}
\os c_4\cd \a_2(c_1)_6\bs(\a_2(c_1)_5)\a_2(c_3)
\a_2(c_1)_4\\  
&=& \a_1(c_1)\os c_2\cd \a_1(c_3)\a_2(c_1)_1\os c_4\cd
\a_2(c_3)\a_2(c_1)_2
\end{eqnarray*}
and
\begin{eqnarray*}
&&\hspace*{-1cm}
 c_{1.-1}c_{2.-1.1}\os c_{1.0}\cd
c_{2.-1.2}\os c_{2.0}\\ 
&=&\sum\a_1(c_{11})\a_1(c_{21})_1\os
c_{12}\cd \a_2(c_{11}) \cd \a_1(c_{21})_2\os c_{22}\cd
\a_2(c_{21})\\ 
&=& \a_1(c_1)\a_1(c_3)_1\os c_2\cd
\a_2(c_1)\a_1(c_3)_2 \os c_4\cd \a_2(c_3)\\ 
&=& \a_1(c_1)\os
c_2\cd \a_1(c_3)\a_2(c_1)_1\os c_4\cd \a_2(c_3)\a_2(c_1)_2
\end{eqnarray*}
\equref{7} follows, and $\tau_\a$ is a  twisting.
\end{proof}

There is also a relation between twisted 2-cocycles and Harrison
2-cocycles. Let $C$ be a right $H$-module coalgebra. Consider
the trivial weak coaction $\rho(c)=1\ot c$, and $\alpha:\ C\to H\ot H$.
The cocycle condition \equref{II} and the twisted comodule condition
\equref{III} of \deref{1.3} then take the following form:
\begin{eqnarray}
&&\hspace*{-2cm}  \a_1(c_2) \os \a_1(c_1)\a_2(c_2)_1\os
\a_2(c_1)\a_2(c_2)_2\eqlabel{d}\\ 
&=& \a_1(c_1)\a_1(c_2)_1\os \a_2(c_1)\a_1(c_2)_2\os
\a_2(c_2)\nonumber\\
&&\hspace*{-2cm} \a_1(c_2)\os \a_2(c_2)\os c_1= \a_1(c_1)\os \a_2(c_1)\os c_2\eqlabel{e}
\end{eqnarray}
The set of Harrison 2-cocyles corresponding to the trivial weak coaction
is denoted by $Z^2_{\rm Harr}(H,C)$. Thus $Z^2_{\rm Harr}(H,C)$ consists of maps satisfying
\equref{I}, \equref{d} and \equref{e}. The set of twisted 2-cocycles 
$\alpha^t:\ C\ot H\to H\ot H$
in the sense of \deref{2.1} will
be denoted by $Z^2_{\rm tw}(H,C\ot H)$.

\begin{proposition}\prlabel{3.4}
Let $C$ be a right $H$-module coalgebra.
We have a bijection between $Z^2_{\rm Harr}(H,C)$ and $Z^2_{\rm tw}(H,C\ot H)$.
\end{proposition}

\begin{proof}
Take $\a^t\in Z^2_{\rm tw}(H,C\ot H)$, and write 
$$\a^t(c\os h)= \sum\a_1^t(c\os h)\os \a_2^t(c\os h)$$
For all $c\in C$ and $h\in H$, we have
\begin{eqnarray}\eqlabel{3.4.1}
&&\a_1^t(c_1\os h_1)\a_1^t(c_2\os
h_2)_1\os
\a_2^t(c_1\os h_1)\a_1^t(c_2\os h_2)_2
\os\\
&&\hspace*{5mm} \a_2^t(c_2\os h_2)=
 \a_1^t(c_1\os h_1)\os \a_1^t(c_2\os h_2)
\a_2^t(c_1\os h_1)_1 \nonumber\\
&&\hspace*{15mm}\os\a_2^t(c_2\os h_2)\a_2^t(c_1\os
h_1)_2\nonumber
\end{eqnarray}
Now define $\alpha:\ C\to H\ot H$ by
$\alpha(c)=\alpha^t(c\ot 1)$.
It is easy to see that $\alpha$ satisfies \equref{I} and \equref{e}.
Using \equref{3.4.1}, we compute
\begin{eqnarray*} 
&&\hspace*{-5mm}  \a_1(c_2)\os \a_1(c_1)\a_2(c_2)_1\os
\a_2(c_1)\a_2(c_2)_2\\ 
&&= \a_1(c_1) \os \a_1(c_2)\a_2(c_1)_1\os
\a_2(c_2)\a_2(c_1)_2\\ 
&&= \a_1^t(c_1\os 1)\os \a_1^t(c_2\os
1)\a_2^t(c_1\os 1)_1\os \a_2^t(c_2\os 1)\a_2^t(c_1\os 1)_2\\ 
&&= \a_1^t(c_1\os 1)\a_1^t(c_2\os 1)_1\os \a_2^t(c_1\os 1)\a_1^t(c_2\os
1)_2\os \a_2^t(c_2\os 1)\\ 
&&= \a_1(c_1)\a_1(c_2)_1\os
\a_2(c_1)\a_1(c_2)_2\os \a_2(c_2)
\end{eqnarray*}
and it follows that $\alpha$ also satisfies \equref{d}.\\
Conversely, let $\a\in Z^2_{\rm Harr}(H,C)$, and define
$\a^t:\ C\ot H\to H\ot H$ by 
$$\a^t(c\os h)= S(h_1)\a_1(c)h_2\os\bs(h_4)\a_2(c)h_3$$
We can easily show that $\alpha^t$ satisfies conditions
\equref{a} and \equref{b} of \deref{2.1}. A straightforward
computation shows that \equref{c} is also satisfied:
\begin{eqnarray*} 
&&\hspace*{-5mm}
 \a_1^t(c_1\os
h_1)\a_1^t(c_3\os h_3)_1\os c_2\os h_2\a_2^t(c_1\os
h_1)\a_1^t(c_3\os h_3)_2\\
&&\hspace*{1cm}\os \a_2^t(c_3\os h_3)\\ 
&=& Sh_1\a_1(c_1)h_2S(h_7)\a_1(c_3)_1h_8\os c_2\os \\
&&\hspace*{1cm}h_5\bs
(h_4)\a_2(c_1)h_3S(h_6)\a_1(c_3)_2h_9\os \bs (h_{11})\a_2(c_3)h_{10}\\
&=& S(h_1)\a_1(c_1)\a_1(c_3)_1h_2\os c_2\os
\a_2(c_1)\a_1(c_3)_2h_3\os\bs (h_5)\a_2(c_3)h_4\\ 
&=& S(h_1)\a_1(c_3)h_2\os c_2\os \a_1(c_1)\a_2(c_3)_1h_3\os\bs
(h_5)\a_2(c_1)\a_2(c_3)_2h_4\\ 
&=& S(h_1)\a_1(c_1)h_2\os c_2\os
\a_1(c_3)\a_2(c_1)_1h_3\os\bs (h_5)\a_2(c_3)\a_2(c_1)_2h_4\\ 
&=& S(h_1)\a_1(c_1)h_2\os c_2\os h_7S(h_8)\a_1(c_3)h_9\bs (h_6)
\a_2(c_1)_1h_3\\
&&\hspace*{1cm}\os\bs (h_{11})\a_2(c_3)h_{10}\bs (h_5)\a_2(c_1)_2h_4\\
&=& \a_1^t(c_1\os h_1)\os c_2\os h_2\a_1^t(c_3\os h_3)
\a_2^t(c_1\os h_1)_1\\
&&\hspace*{1cm}\os \a_2^t(c_3\os h_3)\a_2^t(c_1\os h_1)_2
\end{eqnarray*}
so it follows that $\a^t$ is a twisted 2-cocycle. We leave it
to the reader to show that the maps between $Z^2_{\rm Harr}(H,C)$ 
and $Z^2_{\rm tw}(H,C\ot H)$ defined above are inverses to each other.
\end{proof}

\section{Equivalence of twistings}\selabel{3}
In this Section, we will define an equivalence relation on the set of
twistings of an $H$-module coalgebra $C$. If a twisting is invertible,
then all other twistings in the same equivalence class are also
invertible.

\begin{proposition}\prlabel{2.3}
Take $\tau,\la\in \ct(C)$, and use notation \equref{notation}.
Consider 
$v\in \rh(C,H)$ satisfying the following identities, for all $h\in H,\ c\in C$: 
\begin{eqnarray}
&&
\ve_H\circ v=\ve_C ~~~~;~~~~
v(c\cd h)= S(h_1)v(c)h_2\eqlabel{20}\eqlabel{21}\\ 
&& c_{1.(-1)}v(c_2)_1\os
c_{1.(0)}\cd v(c_2)_2= v(c_1) c_{2.-1}\os c_{2.0.1}\cd
v(c_{2.0.2})\eqlabel{22}
\end{eqnarray}
Then $\psi:\ C^\tau\lra C^\la,~~ \psi(c)= c_1\cd v(c_2)$ is a left
$B$-colinear right $H$-linear coalgebra map inducing the
identity map on $B$. If $v\in {\rm Reg}(C,H)$, then $\psi$ is an
isomorphism.
\end{proposition}

\begin{proof}
Using the second identity in \equref{21} and $B=C/CH^+$, we can easily prove that
$\psi$ is left $B$-colinear and right $H$-linear. Using the first identity
in \equref{20}, we obtain that $\psi$ induces a well-defined map $B\to B$,
which is the identity. In order to prove that $\psi$ is a coalgebra map,
we need to check that
$$ \psi(c_1\cd c_{2.-1})\os
\psi(c_{2.0})= \psi(c)_1\psi(c)_{2.(-1)}\os \psi(c)_{2.(0)}$$ 
Again, we compute the left and right hand side, and see that they are
equal: 
\begin{eqnarray*}
&&\hspace*{-2cm}\psi(c)_1\psi(c)_{2.(-1)}\os
\psi(c)_{2.(0)}\\
&=& (c_1\cd v(c_2))_1(c_1\cd
v(c_2))_{2.(-1)}\os (c_1\cd v(c_2))_{2.(0)}\\ 
&=& c_1\cd
v(c_3)_1(c_2\cd v(c_3)_2)_{(-1)}\os (c_2\cd
v(c_3)_2)_{(0)}\\ 
&=& c_1\cd
v(c_3)_1S(v(c_3)_2)c_{2.(-1)}v(c_3)_3\os c_{2.(0)}\cd
v(c_3)_4\\ 
&=& c_1\cd c_{2.(-1)}v(c_3)_1\os c_{2.(0)}\cd
v(c_3)_2\\ 
&=& c_1\cd v(c_2)c_{3.-1}\os c_{3.0.1}\cd
v(c_{3.0.2})\\ 
&&\hspace*{-2cm} \psi(c_1\cd c_{2.-1})\os \psi(c_{2.0})\\
 &=& (c_1\cd
c_{2.-1})_1\cd v((c_1\cd c_{2.-1})_2)\os (c_{2.0})_1\cd
v((c_{2.0})_2)\\ 
&=& c_1\cd (c_{3.-1})_1v(c_2\cd
(c_{3.-1})_2)\os (c_{3.0})_1\cd v((c_{3.0})_2)\\ 
&=& c_1\cd
c_{3.-1.1}S(c_{3.-1.2})v(c_2) c_{3.-1.3}\os c_{3.0.1}\cd
v(c_{3.0.2})\\ 
&=& c_1\cd v(c_2)c_{3.-1}\os c_{3.0.1}\cd
v(c_{3.0.2})
\end{eqnarray*}
If $v\in {\rm Reg}(C, H)$, then its inverse $w$ also
satisfies \equref{20}, and $\varphi:\ C^{\lambda}\to C^{\tau}$ defined by
$$\varphi(c)= c_1\cd w(c_2)$$
 is the inverse of $\psi$.
\end{proof}

\begin{definition}\delabel{2.4}
We call $\tau,\la \in \ct(C)$ equivalent if
there exists $v\in {\rm Reg}(C,H)$ satisfying
the conditions of \prref{2.3}. We then write $\tau\sim\lambda$.
\end{definition}

\begin{lemma}\lelabel{2.4a}
$\sim $ is  an equivalence relation on $\ct(C)$.
\end{lemma}

\begin{proof}
$\tau\sim \tau$ through $v(c)=\ve(c)1_H$.\\
Next assume that $\tau\sim \la$, and take $v\in {\rm Reg}(C,H)$
satisfying (\ref{eq:20}-\ref{eq:22}). \equref{22} is equivalent to
\begin{equation}\eqlabel{22a}
 c_{(-1)}\os
c_{(0)}= v(c_1) c_{2.-1}v^{-1}(c_3)_1\os c_{2.0.1}\cd
v(c_{2.0.2})v^{-1}(c_3)_2
\end{equation}
The inverse $u$ of $v$ satisfies \equref{20}. It also satisfies
\equref{21} since
\begin{eqnarray*}
&&\hspace*{-2cm}
 u(c_1)c_{2.(-1)}\os c_{2.(0)1}\cd u(c_{2.(0).2})\\ 
&=& u(c_1)v(c_2)c_{3.-1}v^{-1}(c_4)_1\os c_{3.0.1}\cd
v(c_{3.0.3})_1v^{-1}(c_4)_2
\\
&&\hspace*{5mm} S(v^{-1}(c_4)_3)S(v(c_{3.0.3})_2)u(c_{3.0.2})
 v(c_{3.0.3})_3v^{-1}(c_4)_4\\
&=& c_{1.-1}v^{-1}(c_2)_1\os c_{1.0.1}\cd
u(c_{1.0.2})
 v(c_{1.0.3})v^{-1}(c_2)_2\\
&=& c_{1.-1}u(c_2)_1\os c_{1.0.1}\cd u(c_2)_2
\end{eqnarray*}
and it follows that $\la\sim \tau$.\\
Now assume that
$\tau\sim \la$, $\la\sim \ga$, and take the corresponding maps $v,\ u\in {\rm
Reg}(C,H)$.
Set $w=u\ast v$, and write 
$$\tau(c)= c_{-1}\os c_0~~;~~
\la(c)= c_{(-1)}\os c_{(0)}~~;~~\ga(c)= c_{[-1]}\os
c_{[0]}$$
It is easily shown that $w$ satisfies \equref{20}. $v$ satisfies \equref{22a},
and $u$ satisfies
$$  c_{[-1]}\os
c_{[0]}= u(c_1) c_{2.(-1)}u^{-1}(c_3)_1\os c_{2.(0).1}\cd
u(c_{2.(0).2})u^{-1}(c_3)_2$$ 
We compute that
\begin{eqnarray*} 
&&\hspace*{-2cm}
 c_{[-1]}\os
c_{[0]}= u(c_1) v(c_2)
c_{3.-1}v^{-1}(c_4)_1u^{-1}(c_5)_1\os c_{3.0.1}\cd \\ 
&&
v(c_{3.0.3})_1v^{-1}(c_4)_2S(v^{-1}(c_4)_3)S(v(c_{3.0.3})_2)\\
&&\hspace*{5mm}
u(c_{3.0.2}) v(c_{3.0.3})_3v^{-1}(c_4)_4u^{-1}(c_5)_2\\ 
&=& u(c_1)v(c_2)
c_{3.-1}v^{-1}(c_4)_1u^{-1}(c_5)_1\os\\
&&\hspace*{5mm} c_{3.0.1}\cd
u(c_{3.0.2}) v(c_{3.0.3})v^{-1}(c_4)_2u^{-1}(c_5)_2\\ 
&=& (u\ast v)(c_1)c_{2.-1}(u\ast v)^{-1}(c_3)_1\os\\
&&\hspace*{5mm}
c_{2.0.1}\cd (u\ast v)(c_{2.0.2})(u\ast
v)^{-1}(c_3)_2
\end{eqnarray*}
and this proves that $\tau\sim \ga$.
\end{proof}

\begin{theorem}\thlabel{2.5}
Take $\tau\sim\la \in \ct(C)$. If $\tau$ is invertible, then $\la$
is also invertible.
\end{theorem}

\begin{proof}
Take $v\in {\rm Reg}(C,H)$
satisfying the conditions in \prref{2.3}, and let
$\psi:\ C^\tau\lra C^\la$ be the coalgebra isomorphism 
given by
$$\psi(c)= c_1\cd v(c_2)$$
Let $\tau^{-1}$ be the inverse to $\tau$, and write
$$\tau^{-1}(c)= c_{<-1>}\os c_{<0>}~~;~~\tau(c)= c_{-1}\os
c_0~~;~~\la(c)= c_{(-1)}\os c_{(0)}$$
define $\mu:\ C\lra H\os C$ by 
\begin{eqnarray*} 
&&\hspace*{-2cm}
\mu(c)= c_{[-1]}\os c_{[0]}= 
\psi^{-1}(c)_{<-1>}v^{-1}(\psi^{-1}(c)_{<0>1})v(\psi^{-1}
(c)_{<0>3})_1\\
&&\hspace*{1cm}\os \psi^{-1}(c)_{<0>2}\cd v(\psi^{-1}(c)_{<0>3})_2
\end{eqnarray*}
Using the temporary notation $\psi^{-1}(c)_{<-1>}=a$ and $\psi^{-1}(c)_{<0>}=b$, it is not
hard to prove that $\mu$ is a left inverse of $\lambda$.
Indeed, 
\begin{eqnarray*}
&&\hspace*{-2cm}
(\mu\ast\la)(c)=(m\os id)(id\os\la)\mu(c)\\ 
&=& av^{-1}(b_1)v(b_3)_1(b_2\cd v(b_3)_2)_{(-1)} \os (b_2\cd
v(b_3)_{(0)}\\ 
&=& av^{-1}(b_1)v(b_3)_1S(v(b_3)_2)b_{2.(-1)}
 v(b_3)_3\os b_{2.(0)}\cd v(b_3)_4\\
&=& av^{-1}(b_1)b_{2.(-1)}v(b_3)_1\os b_{2.(0)}\cd
v(b_3)_2\\ 
&=& av^{-1}(b_1)v(b_2)b_{3.-1}\os b_{3.0.1}\cd
v(b_{3.0.2})
\end{eqnarray*}
\begin{eqnarray*} 
&=&
 ab_{-1}\os b_{0.1}\cd v(b_{0.2})\\ 
&=& 1\os
\psi^{-1}(c)_1\cd v(\psi^{-1}(c)_2)\\ 
&=&1\os \psi(\psi^{-1}(c))\\
&=&1\os c=\si(c)
\end{eqnarray*}
The proof of the fact that $\mu$ is also a right inverse of $\lambda$
is much more technical.
From the fact that $v$ is invertible, and using \equref{22}, we
obtain
$$\la(c)= c_{(-1)}\os c_{(0)}= v(c_1)c_{2.-1}v^{-1}(c_3)_1 \os
c_{2.0.1}\cd v(c_{2.0.2})v^{-1}(c_3)_2$$
Now set $\psi^{-1}(c)= c_1\cd v^{-1}(c_2)$. We compute
\begin{eqnarray*} 
&&\hspace*{-1cm}(\la\ast\mu)(c)=(m\os id)(id\os
\mu)\la(c)\\ 
&=& v(c_1)c_{2.-1}v^{-1}(c_3)_1 (c_{2.0.1}\cd
v(c_{2.0.2})v^{-1}(c_3)_2)_{[-1]} \\
&&\hspace*{5mm}\os (c_{2.0.1}\cd
v(c_{2.0.2})v^{-1}(c_3)_2)_{[0]}\\ 
&=& v(c_1)c_{2.-1}v^{-1}(c_3)_1S(v^{-1}(c_3)_2)
S(v(c_{2.0.2})_1)(c_{2.0.1})_{[-1]}\\
&&\hspace*{5mm}v(c_{2.0.2})_2v^{-1}(c_3)_3
 \os (c_{2.0.1})_{[0]}\cd
v(c_{2.0.2})_3v^{-1}(c_3)_4\\  
&=& v(c_1)c_{2.-1}S(v(c_{2.0.2})_1)(c_{2.0.1})_{[-1]}
v(c_{2.0.2})_2v^{-1}(c_3)_1 \\
&&\hspace*{5mm}\os(c_{2.0.1})_{[0]}\cd
v(c_{2.0.2})_3v^{-1}(c_3)_2\\ 
&=& v(c_1)c_{2.-1}S(v(c_{2.0.2})_1)\psi^{-1}(c_{2.0.1})_{<-1>}
v^{-1}(\psi^{-1}(c_{2.0.1})_{<0>.1})\\ 
&&\hspace*{5mm} 
v(\psi^{-1}(c_{2.0.1})_{<0>.3})_1 v(c_{2.0.2})_2v^{-1}(c_3)_1\os\psi^{-1}(c_{2.0.1})_{<0>.2}\\
&&\hspace*{5mm} \cd
v(\psi^{-1}(c_{2.0.1})_{<0>.3})_2 v(c_{2.0.2})_3v^{-1}(c_3)_2\\
&=&  v(c_1)c_{2.-1}S(v(c_{2.0.3})_1)(c_{2.0.1}\cd
v^{-1}(c_{2.0.2}))_{<-1>}\\ 
&&\hspace*{5mm} v^{-1}((c_{2.0.1}\cd
v^{-1}(c_{2.0.2}))_{<0>1}) v((c_{2.0.1}\cd
v^{-1}(c_{2.0.2}))_{<0>3})_1\\
&&\hspace*{5mm}
 v(c_{2.0.3})_2v^{-1}(c_3)_1
\os (c_{2.0.1}\cd v^{-1}(c_{2.0.2}))_{<0>2}\\
&&\hspace*{5mm}\cd
v((c_{2.0.1}\cd v^{-1}(c_{2.0.2}))_{<0>3})_2
v(c_{2.0.3})_3v^{-1}(c_3)_2\\ 
&=& v(c_1)c_{2.-1}S(v(c_{2.0.3})_1)S(v^{-1}(c_{2.0.2})_1)c_{2.0.1.<-1>}
v^{-1}(c_{2.0.2})_2\\ 
&&\hspace*{5mm} v^{-1}((c_{2.0.1.<0>}\cd
v^{-1}(c_{2.0.2})_3)_1) v((c_{2.0.1.<0>}\cd
v^{-1}(c_{2.0.2})_3)_3)_1\\
&&\hspace*{5mm}v(c_{2.0.3})_2v^{-1}(c_3)_1 
\os
(c_{2.0.1.<0>}\cd v^{-1}(c_{2.0.2})_3)_2\cd
v((c_{2.0.1.<0>}\\
&&\hspace*{5mm}\cd
v^{-1}(c_{2.0.2})_3)_3)_2v(c_{2.0.3})_3v^{-1}(c_3)_2\\
&=& v(c_1)c_{2.-1}S(v^{-1}(c_{2.0.2})_1v(c_{2.0.3})_1)c_{2.0.1.<-1>}
v^{-1}(c_{2.0.2})_2\\ 
&&\hspace*{5mm} v^{-1}(c_{2.0.1.<0>.1}\cd
v^{-1}(c_{2.0.2})_3) v(c_{2.0.1.<0>.3}\cd
v^{-1}(c_{2.0.2})_5)_1v(c_{2.0.3})_2\\
&&\hspace*{5mm}v^{-1}(c_3)_1\os
c_{2.0.1.<0>.2}\cd v^{-1}(c_{2.0.2})_4\cd
v(c_{2.0.1.<0>.3}\\
&&\hspace*{5mm}\cd
v^{-1}(c_{2.0.2})_5)_2v(c_{2.0.3})_3v^{-1}(c_3)_2
\end{eqnarray*}
\begin{eqnarray*} 
&=& v(c_1)c_{2.-1}S(v^{-1}(c_{2.0.2})_1v(c_{2.0.3})_1)c_{2.0.1.<-1>}
v^{-1}(c_{2.0.2})_2\\ 
&&\hspace*{5mm}
S(v^{-1}(c_{2.0.2})_3)v^{-1}(c_{2.0.1.<0>.1})v^{-1}(c_{2.0.2})_4
S(v^{-1}(c_{2.0.2})_6)_1\\  
&&\hspace*{5mm}
v(c_{2.0.1.<0>.3})_1(v^{-1}(c_{2.0.2})_7)_1v(c_{2.0.3})_2v^{-1}
(c_3)_1\\ 
&&\hspace*{5mm} \os c_{2.0.1.<0>.2}\cd
v^{-1}(c_{2.0.2})_5S(v^{-1}(c_{2.0.2})_6)_2\\ 
&&\hspace*{5mm}
v(c_{2.0.1.<0>.3})_2(v^{-1}(c_{2.0.2})_7)_2v(c_{2.0.3})_3v^{-1}(c_3)_2\\
&=&  v(c_1)c_{2.-1}S(v^{-1}(c_{2.0.2})_1v(c_{2.0.3})_1)c_{2.0.1.<-1>}
v^{-1}(c_{2.0.2})_2 S(v^{-1}(c_{2.0.2})_3)\\
&&\hspace*{5mm} 
v^{-1}(c_{2.0.1.<0>.1})v^{-1}(c_{2.0.2})_4S(v^{-1}(c_{2.0.2})_7)
v(c_{2.0.1.<0>.3})_1\\
&&\hspace*{5mm}v^{-1}(c_{2.0.2})_8v(c_{2.0.3})_2v^{-1}
(c_3)_1\os c_{2.0.1.<0>.2}\cd
v^{-1}(c_{2.0.2})_5\\
&&\hspace*{5mm}S(v^{-1}(c_{2.0.2})_6)
v(c_{2.0.1.<0>.3})_2v^{-1}(c_{2.0.2})_9v(c_{2.0.3})_3v^{-1}(c_3)_2\\
&=&  v(c_1)c_{2.-1}S(v^{-1}(c_{2.0.2})_1v(c_{2.0.3})_1)c_{2.0.1.<-1>}
v^{-1}(c_{2.0.1.<0>.1})\\
&&\hspace*{5mm}v(c_{2.0.1.<0>.3})_1 
v^{-1}(c_{2.0.2})_2v(c_{2.0.3})_2v^{-1}(c_3)_1 \\
&&\hspace*{5mm}\os
c_{2.0.1.<0>.2}\cd v(c_{2.0.1.<0>.3})_2v^{-1}(c_{2.0.2})_3
v(c_{2.0.3})_3v^{-1}(c_3)_2\\  
&=& v(c_1)c_{2.-1}c_{2.0.<-1>}
v^{-1}(c_{2.0.<0>.1})v(c_{2.0.<0>.3})_1v^{-1}(c_3)_1\\  
&&\hspace*{5mm} \os c_{2.0.<0>.2}\cd
v(c_{2.0.<0>.3})_2v^{-1}(c_3)_2\\  
&=& v(c_1)v^{-1}(c_{2.1})v(c_{2.3})_1v^{-1}(c_3)_1 \os
c_{2.2}\cd v(c_{2.3})_2v^{-1}(c_3)_2\\ 
&=& v(c_1)v^{-1}(c_2)v(c_4)_1v^{-1}(c_5)_1 \os c_3\cd
v(c_4)_2v^{-1}(c_5)_2\\
&=& 1\os c=\si(c).
\end{eqnarray*}
and it follows that $\lambda$ is convolution invertible.
\end{proof}

\begin{theorem}\thlabel{3.5}
Let $C$ be a right $H$-comodule algebra, and consider 
$\tau,\ \la \in \ct(C\os H)$. $\tau$ and $\lambda$ are equivalent
in the sense of \deref{2.4} if and only if there is a
left $C$-colinear, right $H$-linear coalgebra isomorphism between the
crossed coproducts $C\smashco _{\alpha} H,\ \rho$ and
$\crth,\ \rho'$ corresponding to $\tau$ and $\lambda$.
\end{theorem}

\begin{proof}
Write
$$\rho(c)= c_{[-1]}\os c_{[0]}~~;~~\rho'(c)= c_{<-1>}\os c_{<0>}$$
If $\tau\sim\la$, then there exists $v\in{\rm Reg}(C\os H,H)$ satisfying
(\ref{eq:20}-\ref{eq:22}). Define 
$$u:\ C\to H~~;~~u(c)=v^{-1}(c\ot 1)$$
If we can show that $u$ satsifies \equref{18} and \equref{19}, then
one implication is proved, by \leref{1.4}. It follows from \equref{22}
that
\begin{eqnarray} 
&&\hspace*{-2cm} (c_1\os 1)_{(-1)}v(c_2\os 1)_1\os (c_1\os
1)_{(0)}\cd v(c_2\os 1)_2\nonumber\\ 
&=& v(c_1\os 1)(c_2\os 1)_{-1}\os (c_2\os
1)_{0.1}\cd v((c_2\os 1)_{0.2})
\eqlabel{3.5.1}
\end{eqnarray} 
applying $1\os 1\os \ve$ to both sides, we find
\begin{eqnarray*} 
&&\hspace*{-2cm} (c_1\os
1)_{(-1)}v(c_2\os 1)\os (1\os \ve)(c_1\os 1)_{(0)}\\
&=& v(c_1\os 1)(c_2\os 1)_{-1}\os (1\os \ve)(c_2\os 1)_0
\end{eqnarray*}
and using \equref{28}, we obtain
$$ c_{<-1>}\os c_{<0>}= u^{-1}(c_1)c_{2[-1]}u(c_3)\os
c_{2[0]}$$ 
so $u$ satisfies \equref{18}.\\
Applying $1\os \ve \os 1$ to both sides of \equref{3.5.1}, we
find
\begin{eqnarray*} 
&&\hspace*{-2cm}
\a'(c)= v(c_1\os 1)(c_2\os 1)_{-1}v^{-1}(c_3\os 1)_1 \os\\
&&(\ve\os 1)((c_2\os 1)_{0.1}\cd v((c_2\os 1)_{0.2}) v^{-1}(c_3\os
1)_2)
\end{eqnarray*}
It follows from \equref{27} that
$$(c\os 1)_0=c_{1[0]}\os\a_2(c_2)$$
and
\begin{eqnarray*} 
&&\hspace*{-2cm}(\ve\os 1)((c_2\os 1)_{0.1}\cd v((c_2\os
1)_{0.2}) v^{-1}(c_3\os 1)_2)\\ 
&=&(\ve\os
1)(c_{2.[0].1}\os\a_2(c_3)_1v(c_{2.[0].2}\os\a_2(c_3)_2)v^{-1}
(c_4\os 1)_2)\\ 
&=&
\a_2(c_3)_1v(c_{2.[0]}\os\a_2(c_3)_2)v^{-1}(c_4\os 1)_2\\ 
&=&
v(c_{2.[0]}\os 1)\a_2(c_3)v^{-1}(c_4\os 1)_2\\ 
&=& v((1\os
\ve)(c_2\os 1)_0\os 1)(\ve\os 1)(c_3\os 1)_0v^{-1}(c_4\os 1)_2
\end{eqnarray*} 
so 
\begin{eqnarray*} 
\a'(c)&=& v(c_1\os 1)(c_2\os
1)_{-1}(c_3\os 1)_{-1}v^{-1}(c_4\os 1)_1\\ 
&&~~~~\os v((1\os
\ve)(c_2\os 1)_0\os 1)(\ve\os 1)(c_3\os 1)_0v^{-1}(c_4\os
1)_2\\ 
&=& u^{-1}(c_1)c_{2[-1]}\a_1(c_3)u(c_4)_1\os
u^{-1}(c_{2[0]})\a_2(c_{3})u(c_4)_2
\end{eqnarray*}
and \equref{19} follows.\\
Conversely, assume that the two crossed coproducts are isomorphic.
By \leref{1.4}, there exists $u\in{\rm Reg}(C,H)$ satisfying
\equref{18} and \equref{19}. Define
$$v:\ C\os H\lra H~~;~~v(c\os h)= S(h_1)u^{-1}(c)h_2$$
Then 
$$\ve_Hv(c\os h)=\ve(S(h_1)u^{-1}(c)h_2)=\ve(c)\ve(h)$$ 
and
\begin{eqnarray*}
v((c\os h)\cd g)&=&v(c\os hg)= S(hg)_1u^{-1}(c)(hg)_2\\
&=& S(g_1)S(h_1)u^{-1}(c)h_2g_2= S(g_1)v(c\os h)g_2
\end{eqnarray*}
so 
\begin{eqnarray*}
&&\hspace*{-1cm}
\la(c\os h)= (c\os h)_1\os (c\os h)_0\\
&=&S(h_1)c_{1.<-1>}\a'_1(c_2)h_2\os c_{1.<0>}\os \a'_2(c_2)h_3\\
&=& S(h_1)u^{-1}(c_1)c_{2.[-1]}u(c_3)u^{-1}(c_4)c_{5.[-1]}\a_1(c_6)
u(c_7)_1h_2\\ 
&&\hspace*{5mm}\os c_{2.[0]}\os
u^{-1}(c_{5.[0]})\a_2(c_6)u(c_7)_2h_3
\end{eqnarray*}
\begin{eqnarray*} 
&=& S(h_1)u^{-1}(c_1)c_{2.[-1]}\a_1(c_3)u(c_4)_1h_2 \os
c_{2.[0].1}\\
&&\hspace*{5mm}\os u^{-1}(c_{2.[0].2})\a_2(c_3)u(c_4)_2h_3\\
&=& S(h_1)u^{-1}(c_1)h_2S(h_3)c_{2.[-1]}\a_1(c_3)h_4S(h_7)u(c_4)_1h_8\\
&&\os c_{2.[0].1}\os
u^{-1}(c_{2.[0].2})\a_2(c_3)h_5S(h_6)u(c_4)_2h_9\\ 
&=& S(h_1)u^{-1}(c_1)h_2S(h_3)c_{2.[-1]}\a_1(c_3)h_4S(h_9)u(c_4)_1h_{10}\os c_{2.[0].1}\\
&&\hspace*{5mm}\os
\a_2(c_3)_1h_5S(h_6)S(\a_2(c_3)_2)u^{-1}(c_{2.[0].2})
\a_2(c_3)_3h_7S(h_8)u(c_4)_2h_{11}\\ 
&=& v(c_1\os h_1)(c_2\os
h_2)_{-1}v^{-1}(c_3\os h_3)_1 \os(c_2\os h_2)_{0.1}\\ 
&&\hspace*{5mm}\cd
v((c_2\os h_2)_{0.2})v^{-1}(c_3\os h_3)_2
\end{eqnarray*}
This shows that  $\tau\sim \la$.
\end{proof}

\section{Twisting Hopf-Galois coextensions}\selabel{4}
Let $H$ be a Hopf algebra with bijective antipode $S$, and $C$ a
right $H$-module coalgebra. As before, we use the following notation
$$B=C/I~~;~~I=\{ c(h-\varepsilon(h))~|~ h\in H,\ c\in C\}$$
For $\tau\in \ct (C)$, we have that $C^\tau/I^\tau=C/I=B$.\\
Now assume that $C/B$ is an $H$-Galois coextension (see \cite{BH}).
This means that the canonical map
$$\beta:\ C\os H\lra C\Box _B C~~;~~\beta(c\os h)= c_1\os c_2\cd h$$
is a bijection.

\begin{lemma}\lelabel{3.1}
With notation as above, consider the map
$$\beta':\ C\os H\lra C\Box _B C~~;~~\beta(c\os h)= c_1\cd h\os c_2$$
If the antipode $S$ is bijective, then $\beta$ is bijective (resp.
injective, surjective) if and only if $\beta'$ is bijective (resp.
injective, surjective).
\end{lemma}

\begin{proof}
The map
$$\phi:\ C\os H\lra C\os H,~~\phi(c\os h)= c\cd h_1\os S(h_2)$$
is a bijection with inverse
$$\phi^{-1}(c\os h)= c\cd h_2\os \bs h_1$$
The statement then follows from the fact that $\beta'=\beta\circ\phi$.
\end{proof}

\begin{theorem}\thlabel{3.2}
Take $\tau\in U(\ct(C))$. Then $C^\tau/B$ is an
$H$-Galois coextension if and only if $C/B$ is an $H$-Galois
coextension.
\end{theorem}

\begin{proof}
Let $\lambda$ be the inverse of $\tau$. As before, we use 
the notation \equref{notation}.
Let $\beta^\tau$ be the canonical map corresponding to the coextension
$C^\tau/B$, that is,
$$\beta^\tau(c\os h)= c_1\cd c_{2.-1}h\os c_{2.0}$$
Consider the following diagram
\begin{equation}\eqlabel{3.2.1}
\begin{diagram}
C\ot H&\rTo^{\beta}&C\Box_B C\\
\dTo^{f}&&\dTo^{g}\\
C\ot H&\rTo^{\beta^{\tau}}&C\Box_B C
\end{diagram}
\end{equation}
where
$$f(c\os h)= c_0\os \bs (c_{-1})h~~;~~
g(c\os d)= c_0\cd \bs (c_{-1})\os d$$
$f$ and $g$ are bijections, with inverses given by
$$f^{-1}(c\os h)= c_{(0)}\os \bs(c_{(-1)})h~~;~~
 g^{-1}(c\os d)= c_{(0)}\cd \bs(c_{(-1)})\os d$$
We can also compute that
\begin{eqnarray*}
&&\hspace*{-2cm}
\beta^\tau f(c\os h)=\beta^\tau( c_0\os
\bs(c_{-1})h)\\  
&=& c_{0.1}\cd
c_{0.2.-1}\bs(c_{-1})h\os c_{0.2.0}\\ 
&=& c_{1.0}\cd
c_{2.-1.2}\bs(c_{1.-1}c_{2.-1.1})h\os c_{2.0}\\ 
&=& c_{1.0}\cd
c_{2.-1.2}\bs(c_{2.-1.1})\bs(c_{1.-1})h\os c_{2.0}\\ 
&=& c_{1.0}\cd \bs(c_{1.-1})h\os c_2\\
&=& c_{1.0}\cd
h_3\bs(h_2)\bs(c_{1.-1})h_1\os c_2\\ 
&=& c_{1.0}\cd
h_3\bs(S(h_1)c_{1.-1}h_2)\os c_2\\
&=&  (c_1\cd h)_0\cd
\bs((c_1\cd h)_{-1})\os c_2\\ 
&=&g( c_1\cd h\os
c_2)=g\beta(c\os h)
\end{eqnarray*}
This shows that \equref{3.2.1} is commutative, and it follows that $\beta$ is bijective if
and only if $\beta^{\tau}$ is bijective.
\end{proof}

\begin{theorem}\thlabel{3.3}
Let $C/B$ be an $H$-Galois coextension, and take $\tau,\la \in \ct(C)$.
Every left $B$-colinear right $H$-linear coalgebra map
$$\psi:\ C^\tau\lra C^\la$$
is of the form
$$\psi(c)= c_1\cd v(c_2)$$
where $v\in \rh(C,H)$ satisfies the conditions (\ref{eq:20}-\ref{eq:22})
of \prref{2.3}. If $\psi$ is an isomorphism, then
$v\in {\rm Reg}(C,H)$.
\end{theorem}

\begin{proof}
We use the notation \equref{notation}. 
As in \cite{BH}, we consider the map
$$\bar{\tau}=(\ve\os 1)\beta^{-1}:\ \cbc\lra H$$
Write $\bar{\tau}(c\os d)=c\ds d$, and recall that $\bar{\tau}$
has the following properties:
\begin{eqnarray} 
\ve_H(c\ds d)&=&\ve_C(c)\ve_C(d)\eqlabel{23}\\
(c\ds d)h&=&c\ds(d\cd h)\eqlabel{24}\\ 
(c\cd h)\ds d&=&S(h)(c\ds d)\eqlabel{25}\\
 c_1\cd (c_2\ds d)&=&\ve(c)d \eqlabel{26}
\end{eqnarray}
The map
$$v:\ C\to H~~;~~v(c)= c_1\ds\psi(c_2)$$
satisfies the property
$$ c_1\cd v(c_2)= c_1\cd (c_2\ds\psi(c_3))=\psi(c)$$
Since $\psi$ is a coalgebra, 
$$\ve_Hv(c)=\ve_H(c_1\ds\psi(c_2))=\ve(\psi(c))=\ve(c)$$
and it follows that $\ve_H\circ v=\ve_C$.\\
It follows from (\ref{eq:25}-\ref{eq:26}) that
\begin{eqnarray*}
&&\hspace*{-2cm}
v(c\cd h)= (c\cd
h)_1\ds\psi((c\cd h)_2) 
= c_1\cd h_1\ds \psi(c_2\cd
h_2)\\ 
&=& S(h_1)(c_1\ds\psi(c_2))h_2= S(h_1)v(c)h_2
\end{eqnarray*}
$\psi$ is a coalgebra map, so
$$ \psi(c_1\cd
c_{2.-1})\os \psi(c_{2.0})= \psi(c)_1\psi(c)_{2.(-1)}\os
\psi(c)_{2.(0)}$$
and 
\begin{eqnarray*}
&&\hspace*{-2cm}
 c_1\cd v(c_2)c_{3.-1}\os c_{3.0.1}\cd v(c_{3.0.2})\\
&&= (c_1\cd v(c_2))_1(c_1\cd
v(c_2))_{2.(-1)}\os (c_1\cd v(c_2))_{2.(0)}\\ 
&&= c_1\cd
v(c_3)_1(c_2\cd v(c_3)_2)_{(-1)}\os (c_2\cd v(c_3)_2)_{(0)}\\
&& = c_1\cd c_{2.(-1)}v(c_3)_1\os c_{2.(0)}\cd v(c_3)_2
\end{eqnarray*}
which is equivalent to
\begin{eqnarray*}
&&\hspace*{-2cm} c_1\os c_2\cd
v(c_3)c_{4.-1}\os c_{4.0.1}\cd v(c_{4.0.2})\\
&=&c_1\os c_2\cd
c_{3.(-1)}v(c_4)_1\os c_{3.(0)}\cd v(c_4)_2
\end{eqnarray*}
After we apply $\beta^{-1}$ to both sides, we obtain
$$ c_1\os v(c_2)c_{3.-1}\os c_{3.0.1}\cd v(c_{3.0.2}) = c_1\os
c_{2.(-1)}v(c_3)_1\os c_{2.(0)}\cd v(c_3)_2$$
and
\begin{equation} 
 v(c_1)c_{2.-1}\os c_{2.0.1}\cd v(c_{2.0.2}) = c_{1.(-1)}v(c_2)_1\os c_{1.(0)}\cd v(c_2)_2\eqlabel{**}
\end{equation}
If $\psi$ is an isomorphism, then $\psi^{-1}:\ C^\la\lra C^\tau$ is
a left $B$-colinear right $H$-linear coalgebra map. Then we have
a map $w:\ C\lra H$ satisfying (\ref{eq:20}-\ref{eq:22}) such that
$$\psi^{-1}(c)= c_1\cd w(c_2)$$
For all $c\in C$, we have that 
$$c=c_1\cd v(c_2)w(c_3)=c_1\cd w(c_2)v(c_3)$$
Proceeding as in the proof of \equref{**}, we find that $v$ is convolution
invertible.
\end{proof}


\begin{thebibliography}{99}
\bibitem{BCZ}
M. Beattie, C.Y. Chen and J.J. Zhang, Twisted Hopf
comodule algebras, {\sl Comm. Alg.} {\bf 24}(5) (1996), 1759--1775.

\bibitem{BT}
M. Beattie and B. Torrecillas, Twistings and Hopf
Galois extensions, {\sl J. Algebra} {\bf 232}(2) (2000), 673--696.

\bibitem{BCM} 
R. Blattner, M. Cohen and S. Montgomery, Crossed products and Inner actions of
Hopf algebras, {\sl Trans. Amer. Math. Soc.} {\bf 298}
(1986), 671--711.

\bibitem{BM}
R. Blattner and S. Montgomery, Crossed products and Galois extensions
of Hopf algebras, {\sl Pacific J. Math.} {\bf 137} (1989), 37--54.

\bibitem{BH}
T. Brzezi\'nski and P.M. Hajac, Coalgebra extensions and algebra
coextensions of Galois type, {\sl Comm. Algebra} {\bf 27}(3) (1999),
1347--1367.

\bibitem{C}
S. Caenepeel, Harrison cohomology and the group of Galois coobjects,
in ``Alg\`ebre non commutative, groupes
quantiques et invariants (Reims, 1995)", 83--101, {\sl S\'emin. Congr.} {\bf 2}, 
Soc. Math. France, Paris, 1997.

\bibitem{CDMP}
S. Caenepeel, S. D\v asc\v alescu, G. Militaru and F. Panaite,
Coalgebra deformations of bialgebras by Harrison cocycles, copairings
of Hopf algebras and double crosscoproducts,
{\sl Bull. Belgian Math. Soc. Simon Stevin} {\bf 4} (1997), 647-671.

\bibitem{CS} 
S. Chase and M. E. Sweedler, ``Hopf algebras and Galois theory'',
{\sl Lect. Notes in Math.} {\bf 97}, Springer Verlag, Berlin, 1969.

\bibitem{DMR}
S. D\v asc\v alescu, G. Militaru and \c S. Raianu, Crossed
coproducts and cleft coextensions, {\sl Comm. Algebra} {\bf 24}(4)
(1996), 1229--1243.

\bibitem{DNR}
S. D\v asc\v alescu, C. N\v ast\v asescu and \c S. Raianu,
``Hopf algebras: an Introduction'', {\sl Monographs Textbooks in Pure
Appl. Math.} {\bf 235}
Marcel Dekker, New York, 2001.

\bibitem{DRZ} S. D\v asc\v alescu, \c S. Raianu and Y.H. Zhang, Finite Hopf Galois
coextensions, crossed coproducts and duality, 
{\sl J. Algebra} {\bf 178} (1995), 400--413.

\bibitem{DT}
Y. Doi and M. Takeuchi, Cleft comodule algebras for a bialgebra,
{\sl Comm. Algebra} {\bf 14} (1986), 801--817.

\bibitem{JC} 
T.X. Ju and C.R. Cai, Twisted Hopf Module
Coalgebras, {\sl Comm. Algebra} {\bf 28}(1) (2000), 307--320.

\bibitem{KT}
H.F. Kreimer and M. Takeuchi, Hopf algebras and Galois extensions
of an algebra, {\sl Indiana Univ. Math. J.} {\bf 30} (1981), 675--691.

\bibitem{Mon}
S. Montgomery, ``Hopf algebras and their actions on rings", American
Mathematical
Society, Providence, 1993.

\bibitem{Sw}
M. E. Sweedler, Cohomology of algebras over Hopf Algebras,
{\sl Trans. Amer. Math. Soc.} {\bf 133} (1968), 205--239.

\bibitem{Swe}
M. E. Sweedler, ``Hopf algebras'', Benjamin, New York, 1969.
\end{thebibliography}
\end{document}